# Chapter 9 From Mathematics and Education, to Mathematics Education

Fulvia Furinghetti, José Manuel Matos, and Marta Menghini

**Abstract** This chapter takes a historical view of the development of mathematics education, from its initial status as a business mostly managed by mathematicians to the birth of mathematics education as a scientific field of research. The role of mathematical communication is analyzed through the growth of journals and research conferences. Actions of internationalization and cooperation in facing instructional and educational problems are illustrated with reference to the journal *L'Enseignement Mathématique* and to ICMI. Curricular and methodological reforms in the 20th century which generated changes in school mathematics are considered. Starting from the acknowledgement that research in mathematics education demands more than the traditional focus on discussing curricular options at distinct grade levels, we identified several specialized clusters, debating specific issues related to mathematics education at an international level. We grouped the clusters into three main areas: relationships with psychology, the study of social, cultural and political dimensions, and the relevance of a theory for mathematics education.

## Introduction

In this chapter we consider the evolution of mathematics education from its initial status as an enterprise mostly managed by mathematicians to the birth of mathematics education as a scientific field of research. We start our story in the 19th century when old states acquired a modern organization, new states were created, and systems of education had to be updated or constructed. In this story researchers in many fields (psychology, philosophy, medicine, sociology, linguistic, anthropology, etc.) had a role, but the main players were professional mathematicians and mathematics teachers. The transition was via a lengthy pathway leading to a clarification of relationships between them and consequent autonomy from mathematicians acquired by mathematics educators. This autonomy was officially acknowledged through the new election procedure for the International Commission on Mathematical Instruction (ICMI)—adopted in Santiago de Compostela (August 19–20, 2006) by the General Assembly of the International Mathematical Union (IMU)—that transferred the election of the ICMI Executive Committee from the IMU General Assembly to the ICMI General Assembly (see Hodgson, 2009). To outline steps taken on this pathway we will focus on the following important moments:

- The attainment by mathematics education of an international dimension at the beginning of the 20th century through the journal *L'Enseignement Mathématique* and the International Commission on the

F. Furinghetti (*) Dipartimento di Matematica, Università di Genova, Genoa, Italy e-mail: furinghe@dima.unige.it

J. M. Matos Faculdade de Ciencias e Tecnologia, Universidade Nova de Lisboa, Lisbon, Portugal e-mail: jmm@fct.unl.pt

M. Menghini Dipartimento di Matematica, Università di Roma Sapienza, Rome, Italy e-mail: marta.menghini@uniroma1.it

Teaching of Mathematics;
- Curricular reforms; and
- The autonomous initiatives inaugurated by the new approach to mathematics education that made it an academic discipline with a new field of research.

Since the shaping of the new discipline benefited from the interaction with other domains, we will also outline the most influential of these interactions.

The movement of communication, internationalization, and solidarity that endowed mathematics education with an international dimension at the beginning of the 20th century involved countries from all around the world, but most events and people that contributed to making mathematics education an academic discipline belonged to Europe and North America. For this reason our history is mainly devoted to these two regions. We leave to the other chapters of this *Handbook* (see, for example, Chapter 26, by Singh and Ellerton) discussion of other rivulets along which mathematics education developed before and after its emancipation from mathematicians and how attention shifted to other regions of the world.

The roots of mathematics education date back to the origins of mankind. Ancient civilizations left us documents that evidence an intertwining between the development of mathematical culture and concern about the transmission of this culture (Karp & Schubring, in press; Kilpatrick, 1994). In the immense landscape of social, economic, and political events that accompanied the evolution of mathematics teaching, we put forward two important developments that affected it. First, the invention of printing in the 15th century created the possibility of universal literacy and for mathematical knowledge to be transmitted easily to large numbers of people. Over centuries, this led to the second development, the creation of schools to educate the masses. As a result, mathematics—which was an arcane subject 600 years ago—has become a subject studied by virtually all students in the world.

For many centuries the roles of mathematics teachers and researchers in mathematics were largely overlapping. Slowly, when mathematical topics reached an advanced stage far from the elementary level, this overlapping happened only in the case of university teachers who were carrying out research as part of their profession. In the primary and secondary schools, however, the division between teaching and researching mathematics became evident. Among other things, this led to a diversified production in mathematical literature: on the one hand, textbooks specifically aimed at school teaching were published; on the other hand, there was a production of materials reporting new results from mathematics research. According to Struik (1987), the process of professionalization of researchers in mathematics was strongly accelerated by the stimuli given to scientific research in the years of the Industrial Revolution, which created "new social classes with a new outlook on life, interested in science and in technical education" (p. 141). New democratic ideas generated by the French Revolution "invaded academic life; criticism rose against antiquated forms of thinking; schools and universities had to be reformed and rejuvenated" (p. 142). In the 19th century the mathematicians' "chief occupation no longer consisted in membership in a learned academy; they were usually employed by universities or technical schools and were teachers as well as investigators" (p. 142).

Around the middle of the 19th century the profession of mathematics teacher at the primary or secondary level was assuming a new shape, in connection not only with the modernization of old nations, but also with the emerging of new social pulses which manifested themselves in new associations and trade unions, political and social movements, and solidarity initiatives. The transmission of mathematical knowledge was no longer a private matter left to families or to religious bodies, but became a public business under the responsibility of the state. In the following years the establishment of modern national systems of instruction took place in the new and old countries. In this process the main concern became the development of curricula, the production of suitable textbooks, and problems associated with teacher education and recruitment. Soon the need for reflecting on problems inherent in the whole construction gave impulse to the creation of specific journals and associations. This is the setting in which our story of the transition from "mathematics and education" to "mathematics education" began.

## Mathematical Communication

Mathematicians, like all scientists, have always felt the need to communicate their results. Towards that end, for a long time they mainly used private communication but after the establishment of academies and societies they began to write proceedings and reports. Following changes in the cultural and social milieu provoked by the Industrial and French Revolutions, the means of communicating became modernized and the first journals devoted specifically to mathematics appeared. Initially they were ephemeral or, like the French *Annales de Mathématiques Pures et Appliquées* lasted for a few decades (1810–1832). But soon, important periodicals, some of them still existing, were published—in 1826, for example, the *Journal für die Reine und Angewandte Mathematik*, founded by August Leopold Crelle, first appeared, and in 1836 the *Journal de Mathématiques Pures et Appliquées*, founded by Joseph Liouville, was published. These journals, and others of the same kind available around that time, contained not only original essays, but also mathematical memoirs extracted from eminent works and abstracts of important papers. In this way they contributed to the progress of mathematics by making available new results and important works not easily accessible to all their readers (among them beginning researchers). They were mainly devoted to research and had an international readership.

Around the middle of the 19th century another kind of journal, usually termed an *intermediate journal*, appeared. Between 1877 and 1881, for example, the *Journal de Mathématiques Élémentaires* [*et Spéciales*] (editor Justin Bourget) was published, and in 1882 it became the *Journal de Mathématiques Spéciales* (editor Gaston Albert Gohierre de Longchamps). Some of these intermediate journals were specifically addressed to teachers and students in classes preparing for admission to special schools, and mathematical themes were treated at an intermediate level between secondary and university. Earlier, in Great Britain, periodicals such as *The Ladies' Diary* or *Woman's Almanack* first issued in 1704 and *The Educational Times*, issued from 1847 to 1929, contributed in some way to the growth of mathematical knowledge by publishing mathematical questions addressed to amateurs. *The Educational Times,* which was linked to the College of Preceptors (1849), developed from the Society of Teachers (1846). This society had been established to improve the standards of secondary school teaching and, according to Howson (2010), "initially offered qualifications for pupils and teachers" (p. 43). It was from this journal that *Mathematical Questions with their Solutions. From "The Educational Times"* (editor William John Clarke Miller) originated, and this was published between 1864 and 1918.

In the panorama of journals of diverse nature appearing in the 19th century it is difficult to identify journals that specifically addressed secondary mathematics teaching. Indeed, a sign that attention would be given to secondary-level mathematics was the presence of the word "elementary" in the title—though the meaning of this term differed in different journals. Examples of the genre were three French publications *Journal de Mathématiques Élémentaires* (editor Henri Vuibert, founded in 1876), *Journal de Mathématiques Élémentaires* (editor de Longchamps, founded in 1882), *L'Éducation Mathématique* (editors Jean Griess and Henri Vuibert, founded in 1898), and the Italian publication *Rivista di Matematica Elementare* (editor Giovanni Massa, founded in 1874). Sometimes the founders and editors of these journals were schoolteachers, and indeed most of the contributors to the Italian *Rivista di Matematica Elementare* were secondary schoolteachers. The dates of foundation show that journals related to mathematics teaching at secondary level were born later than research journals. This delay is understandable if one considers that primary and secondary teachers, who were constructing their professionalism and their identity when the establishment of the systems of education in the various countries was taking place, constituted the main readership of this kind of journal.

The creation of journals devoted to mathematics teaching was often linked with associations of mathematics teachers. In some cases the periodicals provided roots for the idea of founding professional associations. For example, in 1915 in the USA the MAA (Mathematical Association of America) assumed responsibility for the *American Mathematical Monthly*, which was aimed at teachers of mathematics and, since 1894, had been published, privately. The Association of Teachers of Mathematics of the Middle States and Maryland began publishing a quarterly journal, *The Mathematics Teacher* in September 1908, which eventually was adopted as the official journal of the NCTM (National Council of Teachers of Mathematics) upon its founding in 1920. In Italy the journal *Periodico di Matematica* was founded in 1886 and became the official organ of the Italian National Association of Mathematics Teachers, Mathesis. In other cases the founding of teacher associations led to the publication of new journals for the purpose of spreading information and ideas. In Germany the Deutscher Verein zur Förderung des mathematischen und naturwissenschaftlichen Unterrichts was founded in 1891 and the journal *Unterrichtsblätter für Mathematik und Naturwissenschaften* followed in 1895. In the UK the Association for the Improvement of Geometrical Teaching (AIGT), founded in 1871, evolved into the Mathematical Association in 1897. The Association continued to publish *The Mathematical Gazette,* which had first appeared in 1894. In France the APMEP (Association des Professeurs de Mathématiques de l'Enseignement Public), begun its activities and the publication of its *Bulletin* in 1910.

Communication through journals devoted to mathematics (sometimes together with other sciences) accompanied the growth of the community of mathematicians and later of mathematics educators. Some of the 182 mathematics periodicals listed in a *Catalogue* prepared by the Mathematical Association (1913) still survive; and many new ones would be created. Some of these publications primarily addressed mathematical research, but others were devoted to mathematics teaching. The number of the latter grew considerably in the 20th century so that there were 253 in a list compiled by Schubring and Richter (1980). Some of these journals are examined in (Hanna, 2003; Hanna & Sidoli, 2002).

As we mentioned above, in the decades on either side of 1900 most important national associations of mathematics teachers were founded. These associations, and their journals, helped to promote communication and to shape mathematics teacher identity. In particular, the role of the associations was crucial in stimulating and guiding reforms which took place during that period. These reforms worked towards updating school mathematics in accordance with the new trends in research and towards making curricula suitable in an age of industrial and technological innovation. As observed by Nabonnand (2007), it is true that the spirit of reforms was often embodied by strong personalities such as John Perry in the UK, Felix Klein in Germany, and Charles Émile Ernest Carlo Bourlet in France, but the programs of reforms were discussed, worked out and spread with the teacher associations as important players.

In the United States of America the American Mathematical Society (AMS) was formed in 1888. AMS always emphasized research (and still does), whereas MAA emphasized teaching (in colleges), and still does. Eliakin Hastings Moore advocated Perry's ideas in his 1902 AMS Presidential address, see (Moore, 1903), and many AMS members were angry with him for doing so. In the UK, the Association for the Improvement of Geometrical Teaching (AIGT) and, later on, the Mathematical Association, were born with the aim of supporting reforms in the geometric syllabus. In Switzerland, new programs centred on the introduction of graphical representation of functions were introduced following the proposal of the mathematics teachers association. In 1906, after a talk delivered by Emanuel Beke at the annual meeting of the society of Hungarian teachers, a Commission charged with studying general reforms and changes in secondary mathematics teaching was instituted by the same society with Beke as its first president. In certain cases an important role of associations was to defend mathematics teaching when it was marginalized. For example, in Italy the mathematics teacher association *Mathesis*, founded in 1895, had the aim of supporting mathematics teaching against a decline which had started in the 1890s. Most of the associations are still alive and in good health; new ones have been founded. Many publish journals, bulletins, and newsletters, as well as organize national meetings and other activities.

Often, both teachers and professional mathematicians participated in these initiatives. There were also initiatives carried out by secondary teachers alone; this happened, for example, in Italy during the initial period of the teacher associations. In other cases, for example in France, academic mathematicians drove these initiatives and led reform movements. The problem of the relationship between the two communities (mathematicians and mathematics teachers) and the need to share responsibility and authority are ever-present in the background of the development of mathematics education to the status of an academic discipline.

**Mathematics Education Unbounded**

The national journals and teacher associations became an important tool for transmitting ideas and information among teachers within many nations, and proved to be of crucial importance in shaping the identity of mathematics teachers. Considering that the themes treated were related to the national systems of education and that the teachers of a country constituted the readership, it is not surprising that most contributors were national and that the actions of teacher associations were mainly confined to dealing with national problems. In the journals devoted to mathematics teaching the contributions by foreign authors were very few and usually translated into the local language. In spite of these national settings, we can identify some common ground in reflections, at the beginning of the 20th century, on the problems of mathematics teaching. Discussions about the organization of curricula were often based on three main themes:

- Relationship between parts of programs;
- Rigor versus intuition; and
- Relationships between mathematics and the other disciplines.

What emerged was the need to go beyond discussions on the reorganization of the curricula. It was recognized that there was a need to consider new methods of teaching that took into account the following:

- "Practical approaches to teaching," based on observation, experiments and laboratories;
- New findings about children's development; and
- A focus on applications.

Due to many common features among mathematics education problems, possible advantages of international cooperation in working towards solutions to the instructional and other educational problems were recognized in many countries. In the following we will describe two main initiatives that strongly contributed to this growing internationalism.

**The Journal** *L'Enseignement Mathématique*

In the second half of 19th century, internationalization was a perennial idea in many aspects of society. Transportation was becoming speedier, and technological developments facilitated long-distance communication. In this context it was not surprising to see the emergence of the idea of world exhibitions, or fairs, which provided occasions for showcasing new industrial and technological productions and sharing ideas and projects. The first world exhibition was held in London (1851), and this was followed up over the next 30 years with exhibitions in Paris, Vienna, Philadelphia and Melbourne. Internationalism invaded all aspects of life, among them mathematics. It was not by chance, then, that in 1893 a congress of mathematicians was held in Chicago, where a world exhibition was being organized. The 1893 congress of mathematicians was the cornerstone in the process of making mathematics unbounded, and heralded a tradition (started in 1897) of organizing International Congresses of Mathematicians (ICMs). One of the promoters of the tradition of having ICMs was the French mathematician Charles-Ange Laisant, who was stimulated both by his cultural view of the nature of mathematics and by social ideals of fraternity and solidarity (Furinghetti, & Giacardi, 2008).

Following the Congress of Paris, in 1900, ICMs have been held every 4 years (except for breaks due to the two World Wars). These regular forums have contributed remarkably to shaping the identity of an international community of research mathematicians. The International Mathematical Union (IMU) was founded in 1920, and although it was dissolved in 1932 it was re-established in 1951, with the first General Assembly of the new IMU being held in 1952.

The idea of internationalism was not easily transferable into the world of education for two obvious reasons: (a) issues of instruction are mainly national; and (b) mathematics teachers have a status different from that of mathematicians—in particular, they have less opportunities and financial resources for communicating and traveling together. Still, mathematics education was touched by internationalization, thanks to the foundation in 1899 of the journal *L'Enseignement Mathématique* by Laisant and the Swiss mathematician Henri Fehr. The mission and vision of this publication, explicitly declared by the editors in the first issue, was to make mathematics instruction join the movement of solidarity, internationalism and communication of the times.

This international character of *L'Enseignement Mathématique* marked the difference between this journal and the other existing journals addressed to mathematics teaching: immediately, it published surveys on the situation of mathematical instruction in different countries. The editorial board included mathematicians and historians of mathematics who had already shown a genuine interest for the problems of mathematics teaching (notably Klein), and of communication in mathematics (notably Magnus Gustaf Mittag-Leffler, founder of the mathematical journal *Acta Mathematica*).

The early years and the development of *L'Enseignement Mathématique* have been outlined by Furinghetti (2003, 2009). The journal was special not only for its international character, but also for its scope. In the sixth volume (1904) the editors claimed that for them the word "enseignement" (teaching) had the widest possible meaning: it meant teaching to pupils, as well as teaching to teachers—and, indeed, the editors made clear, one can hardly have the one without the other. For this reason they explicitly stated their intention to dedicate a wide place to questions of philosophy, methodology, and history. For them, teachers needed to enlarge their horizons beyond the program of their classrooms and their countries.

*L'Enseignement Mathématique* was a product of the mathematical milieu—but Fehr was teaching in Geneva, where the psychologists Édouard Claparède and Théodore Flournoy were working. They used the journal to launch a questionnaire investigating the ways of working of mathematicians. This study is important because it pointed to aspects that were not merely cognitive—using terminology that we would now say was concerned with the affective domain. On the other hand, research mathematicians, like Henri Poincaré, published articles in the journal that focussed on aspects related to the nature of the mathematical invention.

**The Rise and Development of an International Project: The Early ICMI**

In 1905 David Eugene Smith published in *L'Enseignement Mathématique* a paper that advocated more international cooperation and the creation of a commission to be appointed during an international conference with the aim of studying instructional problems in different countries (see Smith, 1905). This article was the seed for the establishment, during the fourth ICM (Rome, 1908), of the International Commission on the Teaching of Mathematics, with Klein as its first president. In the first decades of its life the Commission was most commonly referred to as CIEM (*Commission Internationale de l'Enseignement Mathématique*), in French, or IMUK (*Internationale Mathematische Unterrichtskommission*), in German. Though it underwent many changes in status and scope, this Commission may be considered the first incarnation of the present ICMI.

The significance of the foundation of ICMI goes beyond the mere creation of an organizational structure. What was important was that it pointed to the existence of an international community for whom the main focus of attention would be mathematics education. Given that the initial members of ICMI were nations, and that the representatives of those nations were predominantly academic mathematicians, it was not surprising that for a long time ICMI's activities were developed inside the community of mathematicians. During ICM meetings, ICMI presented its reports and received mandates for future activities (Furinghetti, 2007; Furinghetti & Giacardi, 2008; Menghini, Furinghetti, Giacardi, & Arzarello, 2008).

The main ICMI outcomes in the early years were national reports on mathematical instruction in the various countries, and international inquiries on important themes of the teaching of mathematics. Although Klein (1923) explicitly claimed that ICMI recognized that all levels of school mathematics deserved attention, in practice attention was mainly paid to secondary and tertiary levels, and to teacher education. These priorities were evident in the following list of activities launched by ICMI between 1908 and 1915:

- Current situation of the organization and of the methods of mathematical instruction;
- Modern trends in the teaching of mathematics; • Rigor in middle school teaching and the fusion of the various branches of mathematics;
- The teaching of mathematics to students of physical and natural sciences;
- The mathematical training of the physicists in the university; •Intuition and experiment in mathematical teaching in the secondary schools;
- Results obtained on the introduction of differential and integral calculus into the upper years of middle school;
- The place and role of mathematics in higher technical instruction; and
- Inquiry into the training of teachers of mathematics in secondary schools in the various countries.

Like many other scientific institutions, ICMI suffered a general crisis during the First World War, and the period between the two world wars was a time of stagnation in ICMI's activities (Schubring, 2008). During the first General Assembly of the reconstituted IMU, held in Rome in 1952, ICMI became a permanent sub-commission of IMU.

However, times had changed, and in the 1950s and 1960s the old agenda based on inquiries and national reports was felt to be inadequate to face new situations. Also, relationships with mathematicians needed to be reconsidered in order to deal with educational problems efficiently.

## Curricular Reforms in the 20th Century

### Reforms at the Beginning of the 20th Century

At the time ICMI was born issues associated with the construction of mathematics curricula were often hotly debated in many countries. These debates not only discussed issues common to the different nations, but also nation-specific matters.

For instance in both the UK and Italy, the adequacy of Euclid's *Elements* for the teaching of geometry was a much-debated topic. In the early 1870s the AIGT (Association for the Improvement of Geometrical Teaching) had been created to consider, and to challenge, the tradition of using rote exercises for the entrance examinations to British Universities. The ensuing discussions generated numerous alternative textbooks and also led to some changes in the entrance examinations. But in 1901 the British Association for the Advancement of Science hosted an address by John Perry that would influence mathematics education throughout the world. Perry attacked the whole system of a mathematical education which, he claimed, did not take into account children's minds, their interests, the applications of mathematics, and connections between different areas of mathematics. His idea of *practical mathematics* applied to the study of geometry meant that the first work with geometry should involve students using rulers, compasses, protractors, set squares, and scissors. In England, the "Perry movement" initiated much discussion about mathematics syllabi and about the need for the reconstitution of secondary mathematics education (Howson, 1982). It also influenced many countries outside England, such as Japan, where Perry had taught for a brief period (Siu, 2009), and the USA, where, as previously mentioned, Moore accepted Perry's arguments and convictions in relation to mathematics education (Moore, 1903).

In Italy an adaptation of Euclid's *Elements* was published in 1867/1868 as the first Italian textbook after the unification. The authors were famous mathematicians who defended the idea of the purity of geometry against criticisms expressed in Italy and in the UK. The Italian reformers emphasized the importance of preparing and publishing good manuals based on the *Euclidean method*. In Italy, research in the field of geometry was flourishing, and many important researchers were engaged in authoring textbooks. For lower secondary school an intuitive geometry was introduced based on observation and on experimental activities.

Towards the end of the 19th century in the USA a "Committee of Ten" was appointed to make recommendations on the standardization, in contents and methods, of American school curricula (Kilpatrick, 1992). The subcommittee for mathematics produced a range of recommendations, for elementary to high school mathematics curricula, which can be summarized in the key words "exercise the pupil's mental activity" and "rules should be derived inductively instead of being stated dogmatically."

In France, a reform of 1902 especially directed at the lycées recognized the need for emphasis on new modern humanities, including mathematics, and to do away with the monopoly of the classical humanities. The reformers also called for school mathematics to take on a greater sense of reality, displaying more applications to the life sciences (Gispert, 2009). An important aspect of the reform was the introduction of elements of differential and integral calculus into secondary schools. We recall that during this period France was a leading country in the field of analysis.

Both the French reform and the Perry movement with its demand for increased emphasis on calculus gave impetus to the German reform movement led by Klein. This movement, whose key phrase was "functional reasoning," had among its principal aims the shifting down of some elements of differential and integral calculus from university to secondary school. However, the contents of the reform were not limited to the last school years; on the contrary, the reform started from the lower grades and involved many teachers (Schubring, 2000). The present-day emphasis given to functions as the conceptual building block for the teaching and learning of algebra and geometry is reminiscent of this German reform movement (Törner & Sriraman, 2005). In particular, the role of analytical geometry in the study of functions was stressed and thus a link between school geometry and algebra was established. Moreover, Klein's Erlanger Program, which characterized geometry as the study of invariant properties under a group of transformations, provided a stimulus for deeper work on geometric transformations in mathematics teaching.

After becoming the foundation president of ICMI in 1908, Klein promoted an international reform based on the ideas of the German reforms. An international comparison of curricula, which was part of ICMI's agenda from the start, was to serve as a key enabling element for this proposal (see Schubring, 2003). Although not all countries participated actively, many initiated significant curriculum reform activities during that period. According to Schubring (2000) these countries included Austria, Belgium, Denmark, France, Germany, Great Britain, Hungary, Sweden, and the USA.

Our analysis of the contents of the mathematics curriculum in various nations led us to concur with Howson (2003) that up to the late 1950s there was considerable agreement on what school algebra might mean. After the introduction of letters to denote numbers or variables should come the construction of algebraic formulae, followed by the formation or solution of linear equations, then quadratics, then simultaneous linear equations, and the properties of the roots of quadratic and cubic equations. In contrast, there might be notable differences in the teaching of geometry. These concerned the closeness to the original Euclid, the level of rigor, the use of algebraic or analytical means, the experimental or intuitive, the use of geometric transformations, and the attention given to space geometry. Nevertheless, it was generally agreed that in most nations attention to a small number of classical theorems in geometry was required—these theorems included the theorem of Pythagoras, the theorem of Thales or intercept theorem, the circle theorems, and congruence and similarity properties.

## Modern/New Math(s)

A second international reform that occurred in the 1960s is thought to have originated from the group of mathematicians established in 1932 under the assumed name Bourbaki. The interest of the Bourbaki group in mathematics education started in the 1950s, when some of its members joined the International Commission CIE A EM

(Commission Internationale pour l' Étude et l' Amélioration de l'Enseignement des Mathématiques), founded by Caleb Gattegno with the aim of studying and improving mathematics teaching (see Félix, 1985). This Commission comprised people from different backgrounds (mathematicians, pedagogists, psychologists, epistemologists, and secondary teachers).

In its initial years, CIEAEM's actions may be summarized in the following points: democratization of mathematics, active pedagogy, and actual involvement of teachers. Among the mathematicians of this research group we find the Bourbakists Jean Dieudonné, Gustave Choquet, and André Lichnerowicz, who also contributed to the text by Piaget et al. (1955), which was the first of the two books edited by CIEAEM. In that book all authors recognized the opportunities that modern mathematics offered in relation to the reform of mathematics teaching, and Dieudonné claimed that the essence of mathematics was reasoning on abstract notions.

The "modern mathematics" movement that developed in Europe had common roots with a parallel movement in the USA (see Moon, 1986)—the new math movement started in the early 1950s by Max Beberman with the creation of the University of Illinois Committee on School Mathematics (UICSM). Soon after the launch of Sputnik in 1957, the American Mathematical Society set up the School Mathematics Study Group (SMSG) to develop a new curriculum for high schools. In 1958, Edward G. Begle, then at Yale University, was appointed as its Director (see Griffiths, & Howson, 1974; Wooton, 1965). Among the many curriculum groups established in the USA during the new math period, SMSG was, perhaps, the most influential. The experiences of this group and the numerous other mathematics curriculum groups established around that time benefited from contributions of psychology (Kilpatrick, 1992).

All these streams of reform related to modern, or new, mathematics met in 1959 at an international conference held in Royaumont, near Paris. The conference was organized by OEEC (Organisation for European Economic Co-operation), and chaired by Marshall Stone, the president of ICMI. An important role was played by members of CIEAEM, particularly by Dieudonné, who gave a lecture concerning the transition from secondary school to university. According to Dieudonné, the treatment of geometry should proceed from the real numbers, establishing rules for the operations on a set of undefined objects so that a vector space structure would be created. Metric relations would then be introduced by means of a scalar product. Euclidean geometry could be dealt with in only three lessons, in which the system of axioms would be presented. The properties of triangles would not have a role in this new development (OEEC, 1961).

We note that in the same year, 1959, the Woods Hole Conference took place in the USA, with the more general aim of improving science education, and bringing together scientists, mathematicians, psychologists and others (Bruner, 1960).

The aim of the Royaumont Conference was to achieve mathematics curriculum reform in Europe—but, since both the USA and Canada had been invited to attend, it could be argued that an international reform stretching beyond European nations was desired. The conference had a more practical sequel in 1962 in Dubrovnik, Yugoslavia, when a group of experts met to produce a modern program for mathematics teaching in secondary schools. In the geometry programs for the ages 15–18 produced by the Commission, the Cartesian plane was defined as a vector space of dimension two with a scalar product. In line with the proposals of Choquet (OEEC, 1962), these concepts were to be introduced via axioms. For children aged from 11 to 15 years, a more intuitive approach to geometry was recommended, in line with the proposals by the Belgian mathematician Paul Libois. So far as algebra was concerned, the contents listed in Dubrovnik included sets, applications and functions, the introduction to real numbers, elements of number theory, combinatorics, groups and structures, linear applications and matrices. Some of these topics would become standard in many curricula. Set theory was to be a major integrating theme, and strongly influenced the language used in textbooks written for modern mathematics.

In both Europe and the USA, the path of innovation was to start at the university and proceed down through the secondary schools to primary schools. Set theory would be present at all levels of education with, for example, cardinal and ordinal aspects of natural numbers being introduced at the beginning of elementary grades (Pellerey, 1989). Many countries officially adopted modern mathematics programs, and in France and Belgium the proposals were completely in line with Bourbakist viewpoints.

Although the modern/new math movements soon aroused strong criticisms (see, e.g., Ahlfors et al., 1962; Kline, 1973; Thom, 1973), the ample debates about changes in school mathematics provided a springboard for subsequent, more solidly based reform initiatives in the 1960s. In the UK the School Mathematics Project was launched in 1961, and the work of Edith Biggs and the "Nuffield Project" popularized the use of concrete materials and of laboratory techniques in British primary school mathematics programs. In 1967 the Nordic Committee for the Modernization of School Mathematics (Denmark, Finland, Norway, and Sweden) presented a new syllabus inspired by new math. One of the best-known members of this Committee was Bent Christiansen, of Denmark. In 1968 the Zentrum für Didaktik der Mathematik (Centre for the Didactics of Mathematics) was founded in Karlsruhe by Hans Georg Steiner and Heinz Kunle. This was followed in 1973 by the IDM (Institut für Didaktik der Mathematik), founded in Bielefeld by Steiner, Michael Otte and Heinrich Bauersfeld, whose aims combined practice in school and theoretical research. In 1969 the first IREMs (Instituts de Recherche sur l'Enseignement des Mathématiques) were established in Lyon, Paris, and Strasbourg. In the early 1970s the Collaborative Group for Research in Mathematics Education was established at the University of Southampton's Centre for Mathematics Education, with Geoffrey Howson and Bryan Thwaites as collaborators. In 1971 Hans Freudenthal founded the Institut Ontwikkeling Wiskunde Onderwijs (IOWO, Institute for the Development of Mathematics Teaching). This initiative had its far roots in the "Mathematics Working Group" founded in 1936 by Tatiana Ehrenfest-Afanassjewa. The meetings of this group were attended by Freudenthal and constituted a first step in the successive development of the "Realistic Mathematics" movement, initially led by Freudenthal (Smid, 2009).

New Bourbakist-type topics such as vectors, transformations, matrices, and set theory were included in the school mathematics curricula of numerous countries, and a greater emphasis on probability and statistics became the order of the day. The 1970s were fertile years for the creation of projects, as shown by the fact that the presentations of 15 projects were mentioned in the *Proceedings of the Third International Congress on Mathematical Education* (ICME-3), held in Karlsruhe, Germany, in 1976. These and other changes in mathematics education were outlined in a special issue of *Educational Studies in Mathematics* entitled "Change in Mathematics Education Since the Late 1950s—Ideas and Realisation: An ICMI Report" (1978).

**Creeping Reforms**

In addition to these strong curricular innovations there were also some creeping reforms that influenced both curriculum content and teaching and learning methods in school mathematics. The experimental work of psychologists, new teaching aids, and the reform movements of the early 20th century brought an interest among mathematicians in mathematics laboratories (Borel, 1904) in which students actively used drawing instruments, calculating machines, and manipulatives. At the beginning of the 20th century, Peter Treutlein, a German mathematician, developed more than 200 models that could assist the teaching of geometry (Treutlein & Wiener, 1912). These models were manufactured and distributed by famous manufacturers such as those of Ludwig Brill (Darmstadt) and Martin Shilling (in Halle and then Leipzig) in the middle of the 20th century, and came to be widely used in German universities and polytechnics.

After the Second World War the use of concrete materials was taken up again in many contexts. In 1945 an NCTM yearbook was devoted to measuring and drawing instruments and to the creation of three-dimensional physical models. An active promoter in this field was Gattegno, who focussed the early activities of CIEAEM on concrete materials (see Gattegno et al., 1958). This activity had an important didactical transposition in the work of the teacher Emma Castelnuovo. Gattegno, as well as the mathematician and psychologist Zoltan Dienes, strongly supported the use of manipulatives, such as Cuisenaire rods and logic blocks, in classroom activities. The presence of Dienes at ICME-1 in Lyon, France, in 1969 testified to the interest of the ICMI community in the use of concrete materials.

Other psychologists, including Jean Piaget, influenced the movement. Willmore (1972) and Price (1995) have pointed out the importance of this in changing thinking about the teaching and learning of mathematics. Libois used concrete materials at the École Decroly in Brussels, and in the UK, the Association of Mathematics Teachers (ATM) strongly supported Gattegno's initiative in promoting the use of manipulatives. Manipulatives became a vehicle for intuition and experiment in the classroom, and prepared school milieu to receive subsequent innovations with mathematical technology (Ruthven, 2008). Gattegno authored innovative software for teaching elementary numeration concepts and films for teaching geometry that extended some of the themes in Jean Nicolet's films (Powell, 2007).

In the *Proceedings* of the first ICME Congress (1969) we find reference to games, worksheets, films, overhead projectors, and to concrete materials to be used in the classroom. The use of materials is put in relation to a new methodology of classroom activities that also includes working groups and classroom discussion. At that time, computers were entering into discussions on mathematics education. An explicit reference to the role of computers in school mathematics, especially for applied mathematics, was made by Bryan Thwaites (1969) in his address at ICME1. At the same conference, Frédérique Papy presented the "minicomputer" (Papy, 1969). The initial interest in the algorithmic aspects or in discrete mathematics created a place for programming to be considered as a means for attaining rigor (Furinghetti, Menghini, Arzarello, & Giacardi, 2008).

In the 1970s and 1980s attention turned towards learning environments, or microworlds, for example, in the form of turtle geometry as presented by Seymour Papert at ICME-2 (Howson, 1973; Papert, 1972a, 1972b). Software was developed, including forerunners to the dynamic geometry software which helped in revitalizing parts of mathematics, for example, proofs and Euclidean geometry. Technology was considered as a means for changing both the curriculum and teaching practices; mathematical activity could be enriched by modelling or processing data in statistics, by experimenting, and by visualizing. Research on the role and use of technology in the teaching of geometry was conducted, at first using a constructivist perspective in a broad sense, and later using additional theoretical perspectives, in particular, the social interactions in which learning takes place (Laborde, 2008). The use of dynamic geometry software was explored as a mediator between constructivist and other theoretical levels, highlighting the need for precise curricular construction (Borba & Bartolini, 2008).

The increasing availability of ordinary calculators, scientific calculators, and graphics calculators generated interesting experimental approaches to instruction. On the one hand, attention was directed at algorithmic aspects (see Engel, 1977), but on the other hand the ways in which some topics—functions, for example—might be dealt with in secondary schools using the new technology began to be investigated (Guin, Ruthven, & Trouche, 2005). Already, at ICME-2 in Exeter, a Working Group had explicitly addressed technology, and at ICME-3 in Karlsruhe this happened with five official activities. The survey presented by Fey (1989) at ICME-6 in Budapest described developments in the use of technology during this pioneering period. The first ICMI Study, launched in 1984, was devoted to computers and informatics (Churchhouse et al., 1986).

## From Mathematics and Education to Mathematics Education

### Emergence of New Approaches in Mathematics Education

In the 1950s mathematical research changed direction, and also the role of mathematics in society changed. New uses of mathematics were promoted by advances in technology, and by the political associations with the space race and the iron curtain. Mathematics instruction was perceived by governments as linked to an important potential for power among nations. In the meantime, schools were being called upon to deal with rapidly increasing populations and associated educational problems.

Given the complexity of emerging educational problems, the mere study and comparison of curricula and programs, which had been the main activities of early ICMI, were judged to be insufficient. New approaches to mathematics education suitable to the changed mathematical and social contexts were needed (Furinghetti, Menghini, Arzarello, & Giacardi 2008). Various initiatives, such as CIEAEM and the USA curricular groups, pointed to the need for cooperation among mathematicians, teachers psychologists, mathematics teacher educators and mathematics teachers. Clearly, there had emerged a need for new professional expertise featuring what Krygowska (1968) called "frontier research," which acknowledged mathematics education as a scientific discipline.

Freudenthal (1963) observed that history had shown the sterility of the problems of mere organization. By the end of the 1960s research interest shifted from curricular issues to the wider study of various dimensions of mathematics education. There emerged a trend towards widening the scope of curricular interventions, for example to pre-school and to vocational and adult education settings. There was also a call for more careful scientific research in mathematics education. A strong case for the importance of empirical research was made in the first ICME in 1969 by Begle, then at Stanford University. According to Begle (1969):

> ... the factual aspect has been badly neglected in all our discussions and ... most of the answers we have been provided have generally had little empirical justification. I doubt if it is the case that many of the answers that we have given to our questions about mathematics education are completely wrong. Rather I believe that these answers were usually far too simplistic and that the mathematical behaviours and accomplishments of real students are far more complex than the answers would have us believe. (p. 233)
>
> Interest in empirical research in mathematics education was growing in the USA and, by the mid-1960s, several conferences discussing priorities for research for mathematics education took place.

By 1968 a Special Interest Group on mathematics education research had been formed within the American Educational Research Association (Kilpatrick, 1992). Although this kind of research was not being embraced in many other countries, the growth of international research journals and centres would change this perspective. As Fehr and Glaymann (1972) stated in the UNESCO publication *New Trends in Mathematics Teaching*:

> The curriculum reform movement of the last two decades in school mathematics was aimed primarily at improving educational practice. It was not designed to increase the number or the quality of research studies in mathematics education. Nevertheless, the reform movement did enormously stimulate such research—in part because curriculum reformers have been asked to demonstrate that their work can make a difference in the classroom; in part because these reformers have recognized that future changes can be managed better if we understand more about the teaching and learning of mathematics; and in part because the ferment in the curriculum has attracted many new scholars to the study of problems in mathematics education. (p. 127)

There was also a growing recognition of the need for the academic legitimacy of specialists in mathematics education to be recognized and respected. The "Resolutions of the First International Congress on Mathematical Education" (1969) assumed that mathematics education was becoming a science in its own right, with its own problems relating to both mathematical and pedagogical content. ICME called for the new science of mathematics education to be given a place in suitable mathematical departments of universities or research institutes.

This discussion about the identity of *mathematics education*, or *didactics of mathematics*—the preferred nomenclature in some countries—was continued at ICME-2 in 1972. Anna Zofia Krygowska, for example, in her contribution to the Working Group on teacher training for prospective secondary teachers, which was chaired by Steiner, identified four aspects of didactics of mathematics: a synthesis of the appropriate mathematical, educational, cultural and environmental ideas; an introduction to research; the nature and situation of the child; and practical experience (see Howson, 1973). Bent Christiansen (1975) distinguished between mathematics education as a process of interaction between teachers and learners in their classes and the didactics of mathematics, which was the study of this process. He recognized in didactics of mathematics the status of a new discipline and pointed out that it must be taught by specialists—"didacticians of mathematics"—and not by general education specialists.

### New Initiatives in Mathematics Education

The rethinking on the role and the methods of mathematics education carried out in the 1950s and the 1960s led to a global discussion that included rethinking about the relationship between mathematicians and mathematics educators and a plan for new ways of communicating among mathematics educators. Two ICMI presidents faced these issues with particular energy—Heinrich Behnke and Freudenthal (Furinghetti & Giacardi, 2010). The former tried to settle administrative relationships, including financial issues, with mathematicians after the rebirth of ICMI in the 1950s and looked for new terms of references. But this was not enough: a cultural cut with mathematicians was necessary and this was made by Freudenthal who acted on the two main issues that were characterizing the dependence on the mathematical community, journals and conferences. Both the initiatives he took were taken independently from IMU.

*L'Enseignement Mathématique*, the official organ of ICMI since its foundation, was becoming a mathematical journal with little room for educational issues. On the other hand, the professional mathematics teaching journals were local and, due to their mission and vision, not suitable for publishing articles on didactic research. So in 1968 Freudenthal founded *Educational Studies in Mathematics* (ESM) (Furinghetti, 2008). According to Hanna (2003), this initiative stimulated other groups to publish mathematics education research journals: *Zentralblatt für Didaktik der Mathematik* (ZDM) (now *The International Journal on Mathematics Education*) was first published in 1969 (the first editors were Emmanuel Röhrl and Steiner) and the *Journal for Research in Mathematics Education* (JRME) was first published in 1970 (the first editor was David C. Johnson). ESM and its contemporary journals would become important vehicles in which questions, methods, and research within the discipline "mathematics education" were developed, reported and discussed.

The development of mathematics education was also accelerated by new ways of meeting at the international level. At its inception ICMI promoted important conferences, such as those in Milan (1911) and Paris (1914) but, because conferences were no longer held during and after World War I, the only scheduled places for discussing didactical issues were in those sections within the quadrennial ICMs that were devoted to the didactics of mathematics. These sections usually encompassed also philosophy of mathematics, history, and logic, and were variously put together or separated according to the inclinations of the organizers. No plenary talk was ever devoted to mathematics education.

In the 1960s the new math movement stimulated some important meetings in the USA and in Europe that focussed on mathematics education research, and ICMI collaborated with UNESCO in organizing some of these conferences. Occasionally the audience was enlarged to include teachers. Freudenthal succeeded in establishing the tradition of having an international conference—the International Congress on

Mathematical Education (ICME)—with regular dates. The first of these conferences (1969 in Lyon) was organized according to a traditional pattern of presenting a sequence of talks, but already at the second of Exeter (UK), in 1972, working groups were organized, and projects presented with the aim of creating the very place for discussion of ideas. Since ICME-3 (Karlsruhe, 1976), ICME meetings have been held on a quadrennial basis.

New perspectives for looking at mathematics education also emerged from within the body of mathematicians. The concluding sentence of the talk delivered by Hassler Whitney, a mathematician who became president of ICMI in 1979, provided evidence that attention might be shifted to the learner:

> We are too used to thinking of the subject matter, and how children can learn it. We must start with the children, to see what they really are. (Whitney, 1983, p. 296)

At ICME-3, in Karlsruhe, the first affiliated study groups were established—HPM (the International Group on the relations between the History and Pedagogy of Mathematics) and PME (the International Group for the Psychology of Mathematics Education) (see Furinghetti & Giacardi, 2008). With these groups a new period began with regular meetings and proceedings. This marked the evolution of the provision of support for researchers in mathematics education.

## Clusters of Specific Issues in the "Discipline" of Mathematics Education

By the middle of the 1970s new tendencies outlined in the last section were manifesting themselves more clearly at the international level. Understanding that the endeavour of searching for directions for mathematics education required more than merely discussing curricular options at the distinct grade levels, ICMI officials met with staff members of UNESCO at the end of 1974 to prepare the elaboration of the fourth volume in the series of books, *New trends in mathematics teaching*. This was an important step towards deepening discussion of issues that had already been raised. The aim was not only to identify major problems in the field of mathematics education but also to guide and monitor the direction and intensity of changes taking place in that field (Steiner & Christiansen, 1979). A methodology favouring in-depth discussion of the chapters was chosen, leading to broader approaches to the issues of mathematics education (D'Ambrosio, 2007). The results of this careful preparation of the book became visible during the third ICME that took place in Karlsruhe in 1976 and constituted a landmark in the history of mathematics education.

As a consequence of this in-depth approach, the fourth volume of *New Trends* contained chapters dedicated to the discussion of curricular issues at various levels—including adult education, university teaching, and the use of technology. These were discussed at a deeper level than ever before and a critical analysis of curriculum development and issues associated with the evaluation of students, teachers and educational materials was presented. The importance of moving on from curricular issues was noticed and appreciated: "Until recently, both research and development had focussed on only one of two main determinants of the learning process: the pupil or the curriculum. They did not consider the influence of the teacher nor of the general context of instruction" (Bauersfeld, 1979, p. 200). The book also contained a chapter on the professional life of teachers of mathematics and another discussing goals and objectives for mathematical education.

The third ICME, at Karlsruhe, and the publication of the fourth volume of *New Trends* have been acknowledged as the starting point for the formation of several specialized clusters of specific issues related to mathematics education at an international level. We will group them into three areas: (a) relationships with psychology; (b) the study of social, cultural and political dimensions; and (c) the relevance of a theory for mathematics education.

## Psychology and Mathematics Education

Since the late 19th century answers to issues related to mathematics teaching and learning have been sought in fields outside of mathematics. Important contributions came from the merging of competencies within various educational sciences and other disciplines: pedagogy, psychology, philosophy, and medicine. The early works carried out in this field concerned pupils with particular needs, but the methods applied in these cases soon proved to be suitable for dealing with problems associated with the teaching and learning of normal children in the primary school. The mathematical content taken into consideration was mainly concerned with arithmetic, but the use of concrete materials affected also the teaching of geometry.

Educators carried out their work in *practice schools* founded and directed with the purpose of experimenting with new teaching methods. In these schools practice was strongly interwoven with research and two different research streams arose: one was concerned with research on teaching methods, the other with the observation of pupil behaviour. In these developments the roots of theories of learning that are concerned with what goes on in the brain of the learner (such as in Piaget's theory) can be recognized, as can theories of instruction that refer to the behaviours a child should undertake in order to learn (such as in Bruner's theory).

The influence of the work of pedagogists and psychologists in mathematics education probably started at the beginning of the 19th century through the Swiss educator Johann Heinrich Pestalozzi. Pestalozzi influenced the teaching and learning of arithmetic and geometry in primary schools in Europe (de Moor, 1995; Howson, 2010) and in the USA (Cajori, 1890). One of his followers was Friedrich Fröbel, the founder of the German kindergarten organization. Fröbel brought his pupils to learn by means of games and other activities—wooden blocks were used to teach arithmetic and concrete geometrical objects to teach geometry.

Johann Friedrich Herbart was another scholar to influence how mathematics was taught in schools. Around 1900, Herbart's ideas influenced elementary teaching and teacher education in various countries (Howson, 1982). Notwithstanding the stages of instruction that Herbart urged teachers to follow (see Ellerton & Clements, 2005), his views of the relationship between teaching and learning can be regarded as being consistent with what later became known as constructivism. It was largely based also on human and social interactions.

The interest in child education grew particularly in the USA as a result of the writings of John Dewey who, in 1896, founded a laboratory school at the University of Chicago. In 1904 Dewey moved to Columbia University, where he spent the rest of his career. Dewey framed all learning as the result of activity. As for mathematics learning, one of his leading premises was that the notion of quantity is grasped by the child as a result of solving practical problems (Stemhagen, 2008). This idea of *active learning* was also present in the work of Maria Montessori, who created a school for children in Rome, and of Ovide Decroly who created the *École de l'Ermitage* in Brussels. Both were physicians who developed their methods when working initially with children with minor disabilities. Decroly's method was based on observations of the surrounding world, but Montessori developed specific materials (materiale strutturato) that were intended to help children to learn autonomously. After that period many psychological laboratories were established in Europe, often by psychologists such as Alfred Binet—the French psychologist famous for his contributions to intelligence theory and testing—and the Swiss neurologist and child psychologist, Claparède. Children's attempts to learn mathematics were often studied in Binet's and Claparède's laboratories.

In the USA, research in the learning of mathematics was conducted by Edward Lee Thorndike, a behaviourist psychologist who had a strong interest in mathematics learning, and William Brownell, a teacher, psychologist, mathematics educator, and education psychologist. Brownell and Thorndike, although coming from different theoretical positions, were part of a broader movement to create a science of education. In 1922 Thorndike published his *Psychology of Arithmetic*, and soon after that his *Psychology of Algebra* (1923). Both were based on the theory of associations in a "connectivist" perspective, and were intended to support Thorndike's series of school mathematics textbooks. Brownell, following the ideas of his advisor Charles H. Judd, stressed the importance of "meaningful learning" with respect to "rote" methods, in contrast to Thorndike's more behaviourist views (Kilpatrick & Weaver, 1977).

Behavioural psychological theories ("behaviourism"), which had been developed via experiments with animals, were linked to school learning by Burrhus Frederic Skinner during the period 1930–1950 with an emphasis being given for what became known as operant conditioning. Skinner emphasized reinforcement processes, seen as fundamental in the shaping of behaviour. According to the corresponding instructional theory, changes in behaviour could be obtained through programmed instruction (or, later on, through mastery learning and computer-assisted learning). These ideas had a wide application in mathematics instruction (see Skinner, 1954), and in particular on theory supporting the early uses of computers in learning.

The major influence of psychology on mathematics education, however, came from the work of the Swiss psychologist, Jean Piaget. While studying the behaviour of children in a clinical manner and identifying "cognitive stages," Piaget developed methods that permitted broadening the range of mathematical topics in primary school. Piaget's stages were paralleled in the USA by the instructional stages of Jerome Bruner but, as Kilpatrick (1992) put it, only "with the arrival of cognitive psychology in 1950s and 1960s, marked by the availability of Piaget's work in English translation and the reinterpretation of that work by Jerome Bruner, [did] researchers in mathematics education begin to have a more judicious regard for psychological theory and to collaborate more frequently with psychologists" (p. 18).

Although the Russian Lev Semënovic Vygotskij was born in the same year as Piaget, it was not until the 1960s that his ideas began to have an impact on mathematics education. This delay was due to the lack of translations of his works and also to a lack of interest in a social perspective in this field. The introduction of Vygotsky's ideas, especially in relation to the crucial role of social interactions in the advancement of learners through their zone of proximal development (ZPD), would prove to be important. For Vygotsky, all knowledge was socially constructed and internalized by joint processes into which learners brought their personal experiences. It followed that close and supportive relationships played an important role in an individual's knowledge growth. In the perspective of cultural mediation, the world of meaning in the child developed by means of tools (artefacts) and signs. Over the past 25 years Vygotskian theory has been applied extensively in mathematics education, the focus being on the mathematical activities of a group of learners or a dyad rather than the individual (Berger, 2005).

An important contribution to the tie between mathematics education and educational sciences came from scholars—such as Caleb Gattegno, Zoltan Dienes, Richard Skemp, and Efraim Fischbein—whose training was both in mathematics and in educational sciences. The work of Skemp and Fischbein stimulated thinking about the role of psychological factors so far as the teaching and learning of mathematics in the higher grades were concerned. Skemp (1976) distinguished between "instrumental" and "relational understanding": Instrumental understanding is the result of a mechanic learning of rules, theorems and their immediate applications, and relational understanding is the result of a personal engagement of the learner with mathematical objects, situations, problems, ideas. We owe to Fischbein deep work on the interactions between intuition and rigor in mathematics education (Tirosh & Tsamir, 2008). Both Skemp and Fischbein were among the founders of PME. Fischbein was the first president of PME, Skemp the second.

During ICME-1, a round table discussion on the psychological problems of mathematics education was organized under the leadership of Fischbein, who also organized and led a similar discussion group at ICME-2. In the introduction to the *Proceedings* for ICME-2, Howson (1973) stressed the importance that Piagetian psychology had in relation to elementary school mathematics. He also noted that the working group on "The Psychology of Learning Mathematics" was the most attended of all working groups at the Congress. According to Howson (1973), the topic discussed "underpins the whole of mathematics education" (p. 15).

In his 1990 introductory chapter providing a research synthesis for PME of the *ICMI Studies Series*, Fischbein (1990) claimed "the psychological problems of mathematical learning and reasoning are scientifically exciting and at the same time genuinely relevant for mathematics education" (p. 4). This sentence epitomized more than a century of interaction between psychologists and mathematics educators. As a matter of fact, though many domains of knowledge have been linked to mathematics education, such as psychology, philosophy, medicine, sociology, linguistic, and anthropology, the main external conceptual support to the development of mathematics education has come from psychology.

**Social, Cultural and Political Dimensions**

In 1972 the chapter dedicated to research in mathematics education in the third volume of *New Trends in Mathematics Teaching* (Fehr & Glaymann, 1972) proposed three areas for research: curricula, methods and materials; learning and the learner; and teaching and the teacher. Four years later, by the time of ICME-3 in Karlsruhe, the chapter on the same issue in the fourth volume (Bauersfeld, 1979) enlarged the possibilities for research activities by listing five possibly fruitful areas for research: investigations of interactions, studies of real classroom situations, research interests of the teacher, extension of the repertoire of research methods, and a theoretical orientation. Events and perspectives presented at ICME-3 were instrumental in mathematics education adopting more comprehensive perspectives.

This widening of prospective research interests was also accompanied by a broader understanding of the dimensions involved in the place and roles of mathematics and mathematics education in society. By the late 1970s there was a growing interest in the importance of social factors either in discussing the role of mathematics in curricula or in the ways in which social and cultural factors intervened in teaching and learning mathematics. It was increasingly recognized that didactics of mathematics is (or should be) "concerned not only with the process of interaction in the classroom but also with mathematics education as a societal aspect: a process of development imbedded in the process of development of the educational system as a whole" (Christiansen, 1975, p. 28). This viewpoint was also expressed elsewhere (for example, Bishop, 1979). Preparation for ICME-3 brought into focus two early tendencies about this theme. One, championed by Ubiratán D'Ambrosio (1979), reflected on the overall objectives and goals of mathematics education; the other was outlined by Bauersfeld (1979), who advocated, among other things, the importance of the study of interactions in the teaching–learning process.

These two areas, a broad perspective of the cultural and social bases for teaching and learning mathematics and the consequent enlargement of the scope of research, saw significant developments in the 1980s. D'Ambrosio's early elaboration of the goals for mathematics education, produced for ICME-3, evolved into a broader perspective offered at his plenary session at ICME-5 in Adelaide (D'Ambrosio, 1985, 2007) when the concept of ethnomathematics was first presented in a major international event in mathematics education. He suggested that mathematics education should take into account the diversity of cultural attitudes and cultural diversity of distinct "societal groups, with clearly defined cultural roots, modes of production and property, class structure and conflicts, and senses of security and of individual rights" (D'Ambrosio, 1985, p. 5). The consideration of the diverse ways in which mathematics blends in distinct cultures and social milieux, together with a reflection of its consequences for mathematics education, prompted a flurry of investigations, many of them uncovering undervalued mathematical activities in daily practices of social groups and professions. This kind of research stimulated further study and reflection on associated educational practices.

Almost at the same time in Europe two lines of research emerged valuing the social dimensions of teaching and learning. Bauersfeld (1980) published his early work about "hidden" social dimensions in the interactions between teacher and students in the mathematics classroom. And Guy Brousseau (1986), immersed in a French tradition of research, proposed a theory accounting for the transformation (and pitfalls) of scientific mathematics into school mathematical knowledge. Both of these lines saw significant developments in further years.

By the end of the 1980s, research on the influence of social and cultural dimensions on mathematics curricula and mathematics teaching and learning was consistently being reported in mathematics education research publications. In a book published in 1987, Stieg Mellin-Olsen, after discussing the mismatch between the mathematical competencies of students in school and in daily life, argued that mathematics education researchers needed to recognize that political dimensions were inevitably at the centre of mathematics teaching and learning (Mellin-Olsen, 1987).

With the benefit of hindsight it can be seen that 1988 was a key year in the development of mathematics education research. During that year, *Educational Studies in Mathematics* dedicated a special issue to "Socio-cultural studies in mathematics education" (Bishop, 1988a, 1988b); Bishop (1988b) authored a book on the subject; a "Fifth Day Special Programme on Mathematics, Education, and Society," at ICME-6 in Budapest, was devoted to "examining the political dimensions of mathematics education" (Keitel, Damerow, Bishop, & Gerdes, 1989, p. i); and, at a plenary at the twelfth PME conference Terezinha Nunes reported on her team's work detailing the mathematical competencies of illiterate children selling small goods in the streets of Brazilian cities (Carraher, 1988). This *social turn*, as Stephen Lerman (2000) called it, signaled "the emergence into the mathematics education research community of theories that see meaning, thinking, and reasoning as products of social activity" (p. 23).

We can include in this social turn the analysis, from an educational stance, of the role of mathematics and mathematics education in society, echoing D'Ambrosio's early reflections on goals for mathematics education (D'Ambrosio, 1979). In the middle of the 1990s, Ole Skovsmose (1994) discussed the relations between mathematics, society, and citizenship. Acknowledging mathematics power in contemporary societies, he proposed the adoption of a critical stance in mathematics education that allowed for a comprehensive perspective connecting issues of globalization, content, and applications of mathematics, as a basis for actions in society, and for empowerment through mathematical literacy.

## A Concern with Theory

The understanding that mathematics education should look for an adequate place in the academic field was already present at the beginning of the 20th century (Kilpatrick, 1992). One of the resolutions passed at the first ICME (1969), related to the need for a "theory of mathematics education" (p. 416). From the middle of the 1970s, in the wake of ICME-3, this push towards theory development became evident. Steiner, based at the IDM at the University of Bielefeld, led this thrust towards theory development and reflection. He formed an international study group called Theory of Mathematics Education (TME), which held five conferences until 1992, and was a regular special group at international conferences. The debate about the nature, the possibilities, the limits and the legitimacy of mathematics education as a scientific field conducted by the group (Steiner et al., 1984) enlarged earlier discussions (e.g., Begle, 1969; Christiansen, 1975) and involved prominent researchers from several countries. The relationship between mathematics education and other fields of knowledge (psychology, education, sociology, mathematics, etc.), the explanatory power of competing paradigms, the viability of home-grown theories, the relationship between theory and practice, and reflections on curriculum change were among the many contributions of this group. The most tangible productions were two books, one edited by Steiner and Vermandel (1988) on the foundations and methodology of mathematics education and another (Biehler, Scholz, Sträßer, & Winkelmann, 1994) offering a comprehensive survey of how mathematics education was viewed around the world.

Several books (Bishop, Clements, Keitel, Kilpatrick, & Laborde, 1996; Bishop, Clements, Keitel, Kilpatrick, & Leung, 2003; English, 2002; Grouws, 1992; Sierpinska & Kilpatrick, 1998; and this *Third Handbook*, in particular) have made an effort to account for the diversity of mathematics education research. In an attempt to characterize this diversity, Bishop (1992, 1998) drew attention to research traditions that were the "result of upbringing education, cultural background, and research training" (Bishop, 1992, p. 712). In 1992, he applied this construct to the characterization of three different traditions and later he used it as the background for a reflection about the relationship between research and educational practice (1998). One tradition is the pedagogue tradition, which values the role of teachers reflecting on their practice, with experiment and observation being the key components of the research. The empirical-scientist tradition was reflected in Begle's paper at the 1969 ICME-1, and "the key to knowledge, and the research process focusses attention on the methods of obtaining that evidence and of analyzing it, often quantitatively" (Bishop, 1992, p. 712). Thirdly, there is the scholastic-philosopher tradition, based on analysis, rational theorizing, and criticism. The actual teaching reality is an imperfect manifestation of these theoretical proposals.